\documentclass[11pt]{article}
\RequirePackage{amsmath,amsthm,amsfonts,latexsym,amscd,amssymb, makeidx, tocbibind,
graphics, showkeys, amsfonts,amssymb,amsmath,amscd,amsthm, graphicx, eufrak, color}
\usepackage{a4wide}
\begin{document}
\begin{center}
\textbf{\large Semi-orthogonal Parseval wavelets associated to GMRAs   \\   on local fields    of positive characteristics}
\end{center}
\begin{center}
Niraj K. Shukla\footnote{Discipline of Mathematics, Indian Institute of Technology Indore, Indore, India.    o.nirajshukla@gmail.com}, Saurabh Chandra Maury\! \footnote{Department of Basic and Applied Sciences, SRM University, Lucknow,  India.   smaury94@gmail.com} and Shiva Mittal\footnote{Department  of Mathematics, S.P.M. Govt. Degree College, Allahabad,     India.   shivamittal009@gmail.com}
\end{center}

\begin{abstract}
 \noindent In this article we establish theory of semi-orthogonal Parseval wavelets associated to generalized multiresolution analysis (GMRA) for the local field of positive characteristics (LFPC). By employing the properties of translation invariant spaces on the core space of GMRA we obtain a characterization of semi-orthogonal Parseval wavelets in terms of consistency equation for LFPC. As a consequence, we obtain  a characterization of an orthonormal (multi)wavelet to be  associated with an MRA in terms of multiplicity  function as well as dimension function of a (multi)wavelet.  Further, we provide characterizations of Parseval scaling functions, scaling sets and bandlimited wavelets  together with a Shannon type multiwavelet for LFPC.
 


\end{abstract}

\begin{flushleft}
\textbf{1. Introduction}
\end{flushleft}

In the recent years wavelets on local fields of positive characteristics (LFPC) are extensively studied by many authors including Behera and Jahan, Benedetto, Debnath and Shah,   Jiang, Li and Ji,   Vyas    and   first author with respect to  multiresolution analysis (MRA),     frame multiresolution analysis, tight wavelet frame,  low-pass filter, etc.  in the references [4-6, 18, 20,  22-24]   but still more concepts need to be studied for its enhancement. Indeed, the development of    theory of wavelet analysis with respect to   groups other than Euclidean spaces, namely,   $p$-adic groups, Cantor dyadic groups,   Vilenkin groups,   locally compact abelian groups (LCAG), Heisenberg group,    etc., was always an interesting part for researchers  due to its various applications [1, 6, 13, 16, 17, 19].  

Our main goal is to develop the theory of semi-orthogonal Parseval (multi)wavelets associated to generalized multiresolution analysis (GMRA) in the setting of LFPC while a rigorous    study of semi-orthogonal Parseval (multi)wavelets and GMRAs has been done by many authors for the case of Euclidean  spaces [2, 3, 8-10].  The concept of GMRA was introduced by  Baggett,   Medina and   Merrill in [2] for separable Hilbert spaces. They used unitary representation of the group of translations acting on the fundamental subspace $V_0$ (known as, \textit{core space}) that generates, through dilation, the subspaces $\{V_j\}_{j \in \mathbb Z}$ associated with the GMRA.   

As the core space $V_0$ possess the properties of translation invariant (TI) spaces and the concept of TI spaces plays a very important role in the development of theory of  GMRAs for the case of Euclidean  spaces [3, 7, 9, 10],    we adopt the approach of TI spaces to establish the theory of semi-orthogonal Parseval multiwavelets associated with GMRAs in the setting of LFPC and obtain   necessary and sufficient conditions for the existence of a semi-orthogonal Parseval (multi)wavelet    in terms of consistency equation.  Further, we obtain a characterization of an orthonormal (multi)wavelet to be  associated with an MRA in terms of dimension function of a (multi)wavelet in the setting of LFPC.  For this, we provide a brief introduction of TI spaces for the case of LFPC while Currey, Mayeli and Oussa in  [14] for nilpotent Lie groups, and Bownik and Ross in [11] for  LCAG  having a co-compact subgroup developed   the theory of TI spaces.  Also  Bownik [7] and Rzeszotnik [21] studied the same  in the  Euclidean setting.

 Although minimally supported frequency   wavelets (wavelet sets) are not well-localized, and thus not directly useful for applications, they have proved to be an essential tool in developing wavelet theory.  The existence of wavelet sets  for LCAG and related groups, and Heisenberg groups  have been discussed  by Benedetto and Benedetto in [6], and Currey and Mayeli in [13], respectively. Wavelet sets are also studied in the references [2, 3,   8,   15, 21].  We also  characterize   bandlimited Parseval   multiwavelet sets of finite order   which in turn  characterizes all multiwavelet sets for LFPC.  Further, we provide a necessary and sufficient condition of scaling functions associated with Parseval multiwavelets which provides characterization of scaling sets.  In this setting, we obtain a Shannon type multiwavelet associated to MRA.

 The present paper is organized as follows:

	 In Section 2, we  provide brief introduction about the  local field of positive characteristics. For this we refer  a  book   by Taibleson [25].  Section 3, is divided into three subsections. In the first subsection, we discuss about core spaces of semi-orthogonal framelets while second subsection contains a characterization of translation invariant spaces along with some basic notions such as range function, multiplicity function and spectral function  for LFPC. In the last subsection, we discuss  semi-orthogonal Parseval multiwavelets, GMRAs and their main characterization theorem that provide  a connection between  the dimension function and the multiplicity function of a (multi)wavelet. We also   prove  that the wavelet multiplicity function satisfies a consistency equation, and the multiplicity function is equal to one in case of multiwavelets associated to MRAs. Finally, in   Section 4, we provide a characterization of bandlimited   Parseval frame multiwavelets for LFPC, and   necessary and sufficient conditions of  Parseval scaling functions   which generalizes  the  characterization of orthonormal scaling functions  provided by Behera and Jahan [4].

\begin{flushleft}
 \textbf{2. Basic Results and  notations of LFPC}
\end{flushleft}

 Throughout the paper, $K$ denotes a local field. By a local field we mean a field which is locally compact, non-discrete, and totally disconnected.   The set $\mathcal O =\left\{x\in K:\ |x|\leq 1\right\}$ denotes the ring of integers which is a unique maximal compact open subring of $K$, where the absolute value  $|x|$ of $x \in K$ satisfies the properties:
$  |x|=0$ if and only if $x=0$;   $  |xy|=|x||y|$, and   $ \ |x+y|\leq \mbox{max}\left\{{|x|,\; |y|}\right\}$, for all $x,\ y\in K$.  Define $\frak P=\left\{{x\in K:\ |x|<1}\right\}$, which is called the \textit{prime ideal} in $K$. In view of totally disconnectedness of $K$, there exists an element $\mathfrak p$ (known as \textit{prime element}) of $\mathfrak P$ having maximum absolute value and then $\mathfrak {P}=\mathfrak {p}\mathcal O$. It can be easily proved that, $\mathfrak P$ is compact and open.  Therefore, the residue space $\mathcal Q= {\mathcal O}/{\mathfrak {P}}$ is isomorphic to a finite field $GF(q)$, where $q=p^c$ for some prime $p$ and positive  integer $c.$

For a measurable subset $E$ of $K$, let $|E|=\int_{K}\chi_{E}(x)dx,$
where $\chi_{E}$ is the characteristic function of $E$ and $dx$ is the Haar measure for $K^+$  (locally compact additive group of K), so $|\mathcal O|=1$. By decomposing $\mathcal O$ into $q$ cosets of $\mathfrak P$, we have $|\mathfrak{P}|=q^{-1}$ and $|\mathfrak{p}|=q^{-1}$, and hence for $x\in K\backslash \{0\}=:K^*$ (locally compact multiplicative group of K), we have $|x|=q^k$, for some $k \in \mathbb Z.$ Further, notice that $\mathcal O^*:=\mathcal O\backslash \frak P$ is the group of units in $K^*$, and for  $x\neq 0$,   we may write $x=\frak p^k x'$ with $x'\in \mathcal O^*$. In the sequel, we denote $\mathfrak{p}^k \mathcal {O}$ by $\mathfrak P^k$, for each $k \in \mathbb Z$ that is known as \textit{fractional ideal}. Here, for    $x \in \mathfrak{P}^k$, $x$ can be expressed uniquely as $x=\sum^{\infty}_{l=k}c_l\ \frak p^l,\ c_l\in \frak U,$ and $c_{k}\neq 0,$ where  $\frak U=\left\{c_i\right\}_{i=0}^{q-1}$ is  a fixed full set of coset representatives of   $\frak P$ in $\mathcal O$.

Let $\chi$ be a fixed character on $K^+$ that is trivial on $\mathcal {O}$ but is nontrivial on $\frak P^{-1}$, which  can be  found   by starting with nontrivial character and rescaling.  For $y\in K$, we define $\chi_y(x)=\chi(yx),   x\in K$.
  For $f\in L^1(K)$,  the \textit{Fourier transform} of $f$ is the function $\hat{f}$ defined by
$$\hat{f}(\xi)=\int_{K}f(x)\overline{{\chi_\xi(x)}}dx=\int_{K}f(x)\chi(-\xi x)dx,$$
which can be extended for $L^2(K)$. 

Notation $\mathbb N_0:=\mathbb N\cup \left\{0\right\}$. Let $\chi_u$ be any character on $K^+$. Since $\mathcal {O}$ is a subgroup of $K^+$, it follows that the restriction $\chi_{u|_{\mathcal {O}}}$ is a character on $\mathcal {O}$. Since $\chi_{u|_{\mathcal {O}}}$ is a character on $\mathcal {O}$, we have $\chi_u=\chi_v$ if and only if $u-v\in \mathcal {O}$.  Hence, we have the following result [25, Proposition 6.1]:

\noindent\textbf{Theorem 2.1.}  \textit{Let $\mathcal Z:=\left\{u(n)\right\}_{n\in \mathbb N_{0}}$ be a complete list of (distinct) coset representation of $\mathcal {O}$ in $K^+$. Then $\left\{\chi_{u(n)|_{\mathcal O}}\equiv \chi_{u(n)}\right\}_{n\in \mathbb N_0}$ is a list of (distinct) characters on $\mathcal {O}$. Moreover, it is a complete orthonormal system on $\mathcal {O}$.}

Next, we proceed to impose a natural order on $\mathcal Z$ which is used to develop the theory of Fourier series on $L^2 (\mathcal O)$. For this,  we choose a set $\left\{1=\epsilon_0,  \epsilon_i\right\}_{i=1}^{c-1}\subset \mathcal O^*$ such that the vector space $\mathcal Q$ generated by $\left\{1=\epsilon_0,\epsilon_i\right\}_{i=1}^{c-1}$ is isomorphic to the vector space $GF(q)$ over  finite field $GF(p)$ of order $p$  as   $q=p^c$. For $n\in \mathbb N_{0}$ such that $0\leq n<q,$ we write $n=\sum_{k=0}^{c-1} a_{k}p^{k},$ where  $0\leq a_{k}<p$. By noting that $\left\{u(n)\right\}_{n=0}^{q-1}$, a complete set of coset representatives of $\mathcal O$ in $\frak P^{-1}$ with $|u(n)|=q$, for $0< n <q$ and $u(0)=0$, we define  $u(n)=(\sum_{k=0}^{c-1} a_k\epsilon_k) \mathfrak{p}^{-1}.$  Now, for $n\geq 0$, we write $n=\sum_{k=0}^s b_{k}q^{k}$, where $0\leq b_{k}<q$, and define $u(n)= \sum_{k=0}^s u(b_k)\frak {p}^{-k}.$ In general, it is not true that $u(m+n)=u(m)+u(n)$ but $u(rq^k+s)=u(r)\frak p^{-k}+u(s)$,  if $r\geq 0, \ k\geq 0$   and  $0\leq s< q^{k}.$

Now, we sum up above in the following theorem (see,  [25, Proposition 6.6],  [4]):

\noindent\textbf{Theorem 2.2.}  \textit{For $n\in \mathbb N_{0}$, let $u(n)$ be defined as above.}
\begin{itemize}
	\item [(a)] \textit{$u(n)=0$ if and only if $n=0$. If $k\geq 1$, then we have $|u(n)|=q^{k}$ if and only if $q^{k-1}\leq n<q^{k}.$}

	\item [(b)]  $\left\{u(k):\ k\in \mathbb N_{0}\right\}=\left\{-u(k):\ k\in \mathbb N_{0}\right\}$.

	\item [(c)]  \textit{For a fixed $l\in \mathbb N_{0}$, we have $\left\{u(l)+u(k):\ k\in \mathbb N_{0}\right\}=\left\{u(k):\ k\in \mathbb N_{0}\right\}$.}
\end{itemize}

  Following result and definition [18]  will be used in the sequel:

\noindent\textbf{Theorem 2.3.} \textit{For all} $l,\ k\in \mathbb N_{0},\ \chi_{u(k)}(u(l))=1.$

\noindent\textbf{Definition 2.4.} A function $f$ defined on $K$ is said to be  \textit{integral periodic} if $$f(x+u(l))=f(x), \ \  \textrm{for \ all} \ l\in \mathbb N_{0}, x\in K.$$

\begin{flushleft}
  \textbf{3. Multiplicity function associated to Semi-orthogonal Parseval wavelets}
\end{flushleft}

  We begin with  Subsection 3.1 by constructing some examples of Parseval multiwavelets along with Shannon type multiwavelet for LFPC and also,  core spaces associated with  semi-orthogonal framelets. Then we provide a brief introduction about translation invariant spaces in Subsection 3.2. Finally, we find  a connection between semi-orthogonal Parseval multiwavelets and generalized multiresolution analysis  by providing a consistency equation  in terms of  multiplicity function of core space in Subsection 3.3, which is the   main aim of this section.  As a consequence, we obtain    an expected result that the multiplicity function associated to an MRA multiwavelets is equal to 1.

\begin{flushleft}
  \textbf{3.1 Core spaces of Semi-orthogonal framelets} 
\end{flushleft}

In  order to  define the concepts of MRA and wavelets on LFPC $K$, we need   notions of translation and dilation. Since $  \bigcup_{j\in \mathbb Z}\frak {p}^{-j}\mathcal O=K,$ we can regard $\frak p^{-1}$ as the dilation (note that $|\frak {p}^{-1}|=q$) and since $\mathcal Z=\left\{u(n):\ n\in \mathbb N_{0}\right\}$ is a complete list of distinct coset representatives of $\mathcal O$ in $K$, the set $\mathcal Z$ can be treated as the translation set. Note that it follows from Theorem 2.2 that the translation set $\mathcal Z$ is a subgroup of $K^+$ even though it is indexed by $\mathbb N_0$.  We  have  the following definition:

  \noindent\textbf{Definition 3.1.1.} A finite set $\Psi=\left\{\psi_m:\ m=1,\:2,...,L\right\}\subset L^2(K)$ is called a \textit{frame wavelet} (simply, \textit{framelet}) in $L^2(K)$ if the system
	$$
	\mathcal A(\Psi):=\left\{D^j_{\mathfrak p}T_{k}\psi_m:\ 1\leq m\leq L, \ j\in \mathbb Z,\ k\in \mathbb N_{0} \right\}
	$$
	forms a    frame   for $L^2(K)$, that means, for each $f \in L^2(K)$, there are $0< A\leq B < \infty$ such that 
		\begin{align*}
		A\|f\|^2\leq  \sum_{m=1}^L\sum_{(j, k) \in \mathbb Z \times \mathbb N_0}\left|<f, D^j_{\mathfrak p}T_{k}\psi_m>\right|^2 \leq B\|f\|^2,
		\end{align*}
		where the dilation and translation operators are defined as follows:
	$$
	D^j_\frak p f(x)=q^{j / 2}f(\frak{p}^{-j}x), \ \mbox{and}\ T_k f(x)=f(x-u(k)), \ \ \textrm{for}\;\; x \in K,  j \in \mathbb Z, k \in \mathbb N_0.
	$$
For $A=B$, 	$\Psi$ is known as \textit{tight frame wavelet} (simply, \textit{tight framelet}) with constant $A$ while it is known as \textit{Parseval    multiwavelet} of order $L$ for $A=B=1$. In case  of Parseval frame system $\mathcal A(\{\psi\})$ for $L^2(K)$, $\psi$ is known as \textit{Parseval    wavelet}. If the system $\mathcal A(\Psi)$ is an orthonormal basis for $L^2(K)$, $\Psi$ is called an \textit{orthonormal multiwavelet (simply,  multiwavelet) of order $L$}  of $L^2(K)$. Moreover, a framelet $\Psi$ is known as \textit{semi-orthogonal} if $D^j_\frak p W \bot D^{j'}_\frak p W $,  for $j \neq j'$, where $W=\overline{\mbox{span}}\{T_k \psi: k \in \mathbb N_0, \psi \in \Psi\}$. 
	
	 For $f \in L^2(K)$ and $j, k   \geq 0$, we have 	\begin{align*} 	T_k D^j_\frak p f(x)=& D^j_\frak p f(x-u(k))=q^{j/2} f(\mathfrak p^{-j} x-\mathfrak p^{-j} u(k))\\ 	=&q^{j/2} f(\mathfrak p^{-j} x-  u(q^j k))=q^{j/2} T_{q^j k} f (\mathfrak p^{-j} x)\\ 	=& D^j_\frak p T_{{q^j k}} f(x), 	\end{align*} which shows that $T_k D^j_\frak p = D^j_\frak p T_{q^j k}$, for all $j, k \geq 0$. Further, we notice that for $f \in L^2(K)$ and $\xi \in K$, we have
$$
\widehat{\left(D^j_{\mathfrak p}T_k f\right)}(\xi)=q^{-j/2} \chi_{u(k)}(-\mathfrak p^j \xi) \widehat{f} (\mathfrak p^j \xi), \ \ \mbox{for} \ j \in \mathbb Z, \ k \in \mathbb N_0.
$$
The following is a necessary and sufficient condition for the system $\mathcal A(\Psi)$ to be a Parseval frame   for $L^2(K)$ [5].

\noindent\textbf{Theorem 3.1.2.}  \textit{Suppose $\Psi=\left\{\psi_m:\ m=1,\:2,...,L\right\}\subset L^2(K).$ Then the affine system $\mathcal A(\Psi)$ is a Parseval frame   for $L^2(K)$
if and only if for a.e. $\xi,$ the following holds}:
\begin{enumerate}
	\item [(i)]  $\displaystyle\sum_{m=1}^L\sum_{j\in\mathbb Z} \left|\widehat{\psi}_m(\mathfrak p^{-j} \xi)\right|^2=1,$  
		\item [(ii)]  $\displaystyle \sum_{m=1}^L\sum_{j\in\mathbb N_0} \widehat{\psi}_m(\mathfrak p^{-j} \xi)\overline{\widehat{\psi}_m(\mathfrak p^{-j} (\xi+u(s))}=0,$ for   $s \in  \mathbb N_0\backslash q \mathbb N_0.$  
\end{enumerate}

In particular, $\Psi$ is an orthonormal  multiwavelet of order $L$ in $L^2(K)$ if and only if $\|\psi_m\|=1$,  for $1 \leq m \leq  L$,  and the above   conditions (i) and (ii)  hold.

\noindent\textbf{Definition 3.1.3.}
A \textit{multiresolution analysis} (MRA) of $L^2(K)$ is a sequence of closed subspaces $\left\{D_{\frak p}^j(V)\right\}_{j\in \mathbb Z}$ of $L^2(K)$ satisfying the following properties:
\begin{enumerate}
	\item [(M1)] $T_k V= V$,  \ \
{(M2)} $V \subset D_{\frak p}(V)$,  \ \ (M3) $\displaystyle\bigcap_{j\in \mathbb Z}D_{\frak p}^j(V)=\left\{0\right\}$,   \ \
(M4) $\overline{\displaystyle\bigcup_{j\in \mathbb Z}D_{\frak p}^j(V)}=L^2(K)$,
\item[(M5)] there is a  $\varphi\in V$ (known as \textit{orthonormal scaling function}) such that $\{T_k \varphi\}_{k\in \mathbb N_{0}}$ forms an orthonormal basis for $V$. 
\end{enumerate}
The space  $V$ is known as \textit{core space}. If we replace (M5) by (M6) as follows: \\ 
(M6) the system  $\{T_k \varphi\}_{k\in \mathbb N_{0}}$ forms a Parseval frame  for $V$, \\ 
then the sequence $\left\{D_{\frak p}^j(V)\right\}_{j\in \mathbb Z}$ is known as \textit{Parseval multiresolution analysis} (PMRA), and $\varphi$ is known as \textit{Parseval scaling function}.

Next, we provide an example  of orthonormal multiwavelet of order $q-1$ which is of \textit{Shannon type} for $L^2(K)$: 

\noindent\textbf{Example 3.1.4.} Let us consider the ring of integers $\mathcal O$ in $K$, and let the set $\left\{u(n)\right\}_{n=0}^{q-1}$ be  a complete set of distinct coset representatives  of $\mathcal O$ in $\frak P^{-1}$ with $u(0)=0$, and $|u(n)|=q$, for $0< n <q$.  Suppose $\Psi=\{\psi_1, \psi_2, \cdots, \psi_{q-1}\}$, is a collection of functions in $L^2(K)$ such that the Fourier transform 
\begin{align*}
\widehat{\psi}_i(\xi)=\chi_{\mathcal O+u(i)} (\xi)=\begin{cases} 1 & \mbox{if} \ \  \ \xi \in \mathcal O+u(i)\\
0 & \mbox{otherwise},
\end{cases}
\end{align*}for all $i, 1 \leq i \leq q-1$. 
Now, we show that $\Psi$ satisfies all the conditions of Theorem 3.1.2. 
\begin{enumerate}
	\item [(a)] {For each $i, 1 \leq i \leq q-1$, $\|\psi_i\|^2=\|\widehat{\psi}_i\|^2=|\mathcal O+u(i)|=|\mathcal O|=1$.}
	\item [(b)]   Since $\mathcal O$ is an additive subgroup  of $\frak P^{-1}$,   we have  
$$
\left\{\mathcal O+u(0), \mathcal O+u(1), \cdots, \mathcal O+u(q-1)\right\},
$$
 a measurable partition of $\frak P^{-1}$, and hence  the system  $  \{\mathcal O+u(1), \cdots, \mathcal O+u(q-1)\} $ is a measurable partition of the set $\frak P^{-1} \backslash \mathcal O=\frak p^{-1} \mathcal O^*$ as $u(0)=0$. As $\displaystyle \bigcup_{j\in \mathbb Z}\frak {p}^{-j}\mathcal O=K,$ and $\mathcal O \subset \frak p^{-1}\mathcal O=\frak P^{-1} $, we have a measurable partition $\{\frak p^{-j} \mathcal O^*: j \in \mathbb Z\}$   of $K$ and hence $\{\mathfrak p^{-j}(\mathcal O+u(i)): j \in \mathbb Z, 1 \leq i \leq q-1\}$ is a measurable partition of $K$. This shows that $\displaystyle\sum_{i=1}^{q-1}\sum_{j\in\mathbb Z} \left|\widehat{\psi}_i(\mathfrak p^{-j} \xi)\right|^2=1,$ for a.e. $\xi$.   
 
\item [(c)]  The term $ \widehat{\psi}_i(\mathfrak p^{-j} \xi)\overline{\widehat{\psi}_i(\mathfrak p^{-j} (\xi+u(s)))}=\widehat{\psi}_i(\mathfrak p^{-j} \xi)\overline{\widehat{\psi}_i(\mathfrak p^{-j} \xi+u(q^j s)))} $ 
is nonzero only when $\mathfrak p^{-j} \xi$ and $\mathfrak p^{-j} \xi+u(q^j s)$ are in the support  of $\widehat {\psi}_i$, denoted by $supp \ \widehat {\psi}_i$,  for $1 \leq i \leq q-1.$ But, for $s \in  \mathbb N_0\backslash q \mathbb N_0$ and $1 \leq i \leq q-1$, both $\mathfrak p^{-j} \xi$ and $\mathfrak p^{-j} \xi+u(q^j s)$ will not be members of   $supp \ \widehat {\psi}_i  =\mathcal O+u(i)$ since   $\{\mathcal O+u(k): k \in \mathbb N_0 \}$ is a measurable partition of $K$. Hence we have   
$$\displaystyle \sum_{i=1}^{q-1}\sum_{j\in\mathbb N_0} \widehat{\psi}_i(\mathfrak p^{-j} \xi)\overline{\widehat{\psi}_i(\mathfrak p^{-j} (\xi+u(s)))}=0, \ \ \mbox{for} \    s \in  \mathbb N_0\backslash q \mathbb N_0.
$$ 
 This proves that, $\Psi$ is an orthonormal multiwavelet of order $q-1$.

\end{enumerate}

Following are examples of Parseval   multiwavelets of order 1 as well as  $q-1$    in $L^2(K)$:

\noindent\textbf{Example 3.1.5.} \textbf{(a).} Let $m \in \mathbb N$.  Suppose $\psi$ is a  function in $L^2(K)$ whose Fourier transform is defined  as follows:
\begin{align*}
\widehat{\psi} (\xi)=\chi_{\frak p^{m} \mathcal O^*} (\xi)=\begin{cases} 1 & \mbox{if} \ \  \ \xi \in \frak p^{m} \mathcal O^*=\frak P^{m}   \backslash \frak P^{m+1}\\
0 & \mbox{otherwise}.
\end{cases}
\end{align*}
Then, $\psi$ is a Parseval     multiwavelet of order $1$ in view of the following: 
  \begin{enumerate}
	 \item [(i)]   The system $\{\frak p^{j} (\frak p^{m} \mathcal O^*): j \in \mathbb Z\}$ is a measurable partition of $K$ since $\displaystyle \bigcup_{j \in \mathbb Z} \frak p^{-j} \mathcal O=K,$     $\mathcal O \subset \frak P^{-1}$, and $\frak p^{-1} \mathcal O^*=  \frak P^{-1}\backslash \mathcal O$. 
\item [(ii)]  The system $\{ \frak p^{m} \mathcal O^*+ u(k): k \in \mathbb N_0\}$ is a measurable partition of a measurable subset of $K$ since $\{\mathcal O+ u(k): k \in \mathbb N_0\}$ is a   measurable partition of $K$ and $\frak P^m \subset  \mathcal O.$
 \end{enumerate}
 \textbf{(b).} Let $m \in \mathbb N$.  Suppose $\Psi=\{\psi_1, \psi_2, \cdots, \psi_{q-1}\}$ is a collection of functions in $L^2(K)$ whose Fourier transforms are defined as follows for each $i, 1 \leq i \leq q-1$:
\begin{align*}
\widehat{\psi}_i(\xi)=\chi_{\mathfrak {P}^m+\mathfrak {p}^m u(i)} (\xi)=\begin{cases} 1 & \mbox{if} \ \  \ \xi \in \mathfrak {P}^m+\mathfrak {p}^m u(i)\\
0 & \mbox{otherwise}.
\end{cases}
\end{align*}
 Then, in view of Example 3.1.4, with   $\mathfrak {P}^m = \mathfrak {p}^m\mathcal O$,     $\mathfrak {P}^m \subset \mathcal O$, and the system $\{\frak p^m W_i: 1 \leq i \leq q-1\}$ is a measurable partition of $\frak p^{m-1} \mathcal O^*$, where $W_i=\mathcal O + u(i)$, $\Psi$ is a Parseval     multiwavelet of order $q-1$.  Here, note that $|\mathfrak {P}^m|=\frac{1}{q^m} <1,$ as $q\geq 2$, and for $k, k' \in \mathbb N_0,  (k \neq k')$, we have  
\begin{align*}
|(\frak p^{m} W_i+u(k)) \cap (\frak p^{m} W_i+u(k'))|=&q^m |(W_i+\frak p^{-m} u(k)) \cap (W_i+\frak p^{-m} u(k'))|\\
=& q^m |(W_i+ u(q^m k)) \cap (W_i+ u(q^m k'))|\\
=&0,
\end{align*}
as the system $\{W_i+ u(k): k \in \mathbb N_0\}$ is a   measurable partition of $K$. 
 
  Now, we state a lemma which  is a basic perturbation result for frames stated in [12] that will be useful to  construct  examples of  framelets but not Parseval (multi)wavelets   for LFPC with the help of above examples:

   \noindent\textbf{Lemma 3.1.6.} \textit{Suppose that $H$ is a Hilbert space, $\{f_j\} \subset  H$ is a 
   frame with constants $C_1$ and $C_2$,
$$
C_1\|f\| \leq \sum_j |<f, f_j>|^2 \leq C_2 \|f\| \  \  \mbox{for all} \  f \in  H,
$$   
  and $\{g_j\}\subset  H$ is a Bessel sequence with constant $C_0$,
$$
\sum_j |<f, g_j>|^2 \leq C_0   ||f ||_2 \  \  \  \mbox{for all} \  f\in  H.
$$
If $C_0 < C_1$, then $\{f_j + g_j\}$ is a frame with constants $((C_1)^{1/2}- (C_0)^{1/2})^2$  and $((C_2)^{1/2} + (C_0)^{1/2})^2$.}

 In Example 3.1.4, the core space $V$ of Shannon type multiwavelet is given by 
	$$
	V=\overline{\mbox{span}}\{\varphi(\cdot-u(k)): k \in \mathbb N_0, |\widehat{\varphi}|=\chi_{S}\},  
	$$ 
	where $S=\mathcal O$. Notice that the sequence $\{D^j_{\mathfrak p} V\}_{j \in \mathbb Z}$  satisfies all the axioms of MRA.  While the  core spaces $V$ and $V'$ of  part (a) and part (b)  of Example 3.1.5 are defined as follows:
		\begin{align*}
		V=&\overline{\mbox{span}}\left\{\varphi(\cdot-u(k)): k \in \mathbb N_0, |\widehat{\varphi}|=\chi_{\frak P^{m+1}}\right\}, \ \ \mbox{and}\\
	V'=&\overline{\mbox{span}}\left\{\varphi(\cdot-u(k)): k \in \mathbb N_0, |\widehat{\varphi}|=\chi_{\bigcup_{i=1}^{q-1}\bigcup_{j \in \mathbb N}\frak p^{j} (\frak p^{m}W_i)}\right\}, 	
		\end{align*}
respectively in which $\{D^j_{\mathfrak p} V\}_{j \in \mathbb Z}$ and $\{D^j_{\mathfrak p} V'\}_{j \in \mathbb Z}$  satisfy  all the axioms of MRA  other than axiom $(M5)$.

Such sequences $\{D^j_{\mathfrak p} V\}_{j \in \mathbb Z}$ satisfying conditions $(M1)-(M4)$ are known as \textit{generalized Multiresoluton analyses} (GMRA), the concept of which was introduced by Baggett, Medina and Merill in [2] for separable Hilbert spaces. They developed GMRA structure for $L^2(\mathbb R^n)$. Since the core space $V$ possess the property of TI space, it motivates us to use the theory of TI spaces for $L^2(K)$ for developing connection between the GMRA structure and (multi)wavelets or framelets for $L^2(K)$. 
 
Given a finite family $\Psi \subset L^2(K)$, we define  its \textit{space of negative dilates}  $V$ as follows:
$$
V =\overline{\mbox{span}}\{q^{j/2}\psi(\frak p^j\cdot-u(k)): j<0, k\in \mathbb N_0, \psi\in \Psi\}.
$$
    
We say that a framelet  $\Psi$ comes from a GMRA if its space of negative dilates   $V $ satisfies $(M1)-(M4)$. In addition, if $V$ satisfies $(M5)$, then $V$ is associated with an MRA.
 
   Since for a semi-orthogonal framelet  $\Psi$, its space of negative dilates   $V $ and the space $W=\overline{\mbox{span}}\{T_k \psi: k \in \mathbb N_0, \psi \in \Psi\}$ satisfy
$$
\bigoplus_{j \in \mathbb Z} D^j_\frak p W=L^2(K), \ \ \ V=\bigoplus_{j\leq -1} D^j_\frak p W=\left(\bigoplus_{j\geq 0} D^j_\frak p W\right)^{\bot},
$$ 
it can be easily seen that  every semi-orthogonal framelet $\Psi$ comes from a GMRA.

A study of translation invariant spaces will be useful to see its converse, that is, when a GMRA gives rise to a (multi)wavelet, or a semi-orthogonal framelet.

  \begin{flushleft}
    \textbf{3.2.  Translation invariant spaces for LFPC}
     \end{flushleft}

Suppose  $\mathcal Z=\{u(k): k \in \mathbb N_0\}$. A closed subspace $V$ of $L^2(K)$ is said to be \textit{translation invariant} under $\mathcal Z$, in short $\mathcal Z$-TI if $f\in  V$
implies $T_k f \in  V$ for all $k \in \mathbb N_0$. Given a countable  subset $\mathcal A\subset L^2(K)$, we define a $\mathcal Z$-TI space generated by $\mathcal A$ as
$$\mathcal S^{\mathcal Z}(\mathcal A)=\overline{\mbox{span}}\left\{T_k f : f\in \mathcal A, k \in \mathbb N_0\right\}\subset L^2(K).$$
If $\mathcal A=\{\varphi\}$, then the space $\mathcal S^{\mathcal Z}(\{\varphi\})$ is called  \textit{principal translation  invariant} under $\mathcal Z$ ($\mathcal Z$-PTI) which is denoted by $V_{\varphi}$.

 It is easy to see in view of  Zorn's lemma that any $\mathcal Z$-TI space  can be decomposed into an orthogonal sum of $\mathcal Z$-PTI spaces. That means, if   $V$ is a $\mathcal Z$-TI space, then there exists a countable set of functions $\{\varphi_i\}_{i=1}^M$ belonging to $V$ such that
$V=\oplus_{i=1}^M V_{\varphi_i},$ where   $M$ is a natural number or infinity. Here, the orthogonality condition of   $\cal Z$-PTI spaces, namely, $V_{\varphi_1}$ and $V_{\varphi_2}$ can be written as :
$$
\sum_{k \in \mathbb N_0} \widehat{\varphi}_1(\xi+u(k))\overline{\widehat{\varphi}_2(\xi+u(k))}=0, \ a.e.\ \xi \in K,
$$
for  $\varphi_1, \varphi_2 \in L^2(K)$.  This   follows by noting that  $V_{\varphi_{1}}$ is orthogonal to $V_{\varphi_{2}}$ if and only if $ <T_k\varphi_1,\varphi_2> = 0,$ for all $k\in \mathbb N_0.$ Now,  using the argument of standard periodization,   we obtain $$ <T_k\varphi_1,\varphi_2>  
 = \int_{\mathcal O}\chi_{u(k)}(\xi)\sum_{p\in \mathbb N_{0}}\hat{\varphi_1}(\xi+u(p)) \overline{\hat{\varphi_2}(\xi+u(p))} d{\xi}=0,$$ for all $k\in \mathbb N_{0}$, and hence by using  the  property  of Fourier coefficients of the integral periodic function we conclude the condition. 

 Next theorem plays an important role  to find building blocks of all $\mathcal Z$-TI spaces.

\noindent\textbf{Proposition 3.2.1.} \textbf{(A)} \textit{Let $V_{\varphi}$ be a  $\mathcal Z$-PTI space. Then  $f \in V_\varphi$ if and only if $\widehat{f}(\xi)=r(\xi)\widehat{\varphi}(\xi),$  
for some integral periodic function $r \in L^2\left(\mathcal O, w\right)$, where
\begin{align*}
w(\xi)=\sum_{k\in \mathbb N_0}\left|\widehat{\varphi}\left(\xi+u(k)\right)\right|^2.
\end{align*}
The support of $w$ is known as  the  \textit{spectrum} of $V_\varphi$ which is  denoted by $\Omega$.}\\
\textbf{(B)} \textit{Let  $\varphi \in L^2(K)$. Then a necessary and sufficient condition for the system $\left\{T_k\varphi: k\in \mathbb N_0\right\}$ to be a Parseval  frame  for the $\cal Z$-PTI space  $V_\varphi$ is as follows:
\begin{align*}
\sum_{k\in \mathbb N_0}\left|\widehat{\varphi} (\xi+u(k))\right|^2= \chi_{\Omega}(\xi), \ \ a.e. 
\end{align*}
 }
\noindent\textbf{Proof.} \textbf{(A)} This follows by noting  that $V_\varphi=\overline{\mathcal A_\varphi}$, $L^2\left(\mathcal O, w\right)=\overline{\mathcal P_\varphi}$ and the operator $U: \mathcal A_\varphi \rightarrow \mathcal P_\varphi$ defined by $U(f)(\xi)=r(\xi)$ is an isometry which is onto, where  $\mathcal A_\varphi=\mbox{span}\left\{T_k\varphi: k \in \mathbb N_0\right\},$ and  $\mathcal P_\varphi$ is the space of all integral periodic trigonometric polynomials $r$ with the $L^2\left(\mathcal O, w\right)$ norm
$$
\|r\|^2_{L^2\left(\mathcal O, w\right)}=\int_{\mathcal O} |r(\xi)|^2 w(\xi) d\xi, \ \ \mbox{where}  \quad  w(\xi)=\sum_{k\in \mathbb N_0}\left|\widehat{\varphi}\left(\xi+u(k)\right)\right|^2.
$$
Here, $f \in \mathcal A_\varphi$ if and only if  for  $r\in \mathcal P_\varphi$, $\widehat{f}(\xi)=r(\xi)\widehat{\varphi}(\xi),$ where $r(\xi)=\displaystyle\sum_{k\in \mathbb N_0} a_{k} \overline{\chi_{u(k)}(\xi)}$, for  a finite number of non-zero  elements of $\{a_{k}\}_{k \in \mathbb N_0}$.  Now, by  {splitting} the integral into cosets of $\mathcal O$ in $K$ and using the fact of integral periodicity of  $r$, we have
\begin{align*}
\|f\|^2_2=\int_{\mathcal O} \sum_{k\in \mathbb N_0} \left|\widehat{f}\left(\xi+u(k)\right)\right|^2 d\xi
=\int_{\mathcal O} |r(\xi)|^2\sum_{k\in \mathbb N_0} \left|\widehat{\varphi}\left(\xi+u(k)\right)\right|^2 d\xi
=\|r\|^2_{L^2\left(\mathcal O, w\right)},
\end{align*}
 which shows that the operator $U$ is an isometry. \\
\textbf{(B)} Notice that  for every $f \in V_\varphi$, we have $\widehat{f}(\xi)=r(\xi) \widehat{\varphi}(\xi),$
for some integral periodic function $r \in L^2\left(\mathcal O, w\right)$,
and hence
\begin{align*}
\sum_{k \in \mathbb N_0} \left|<f, T_k\varphi>\right|^2
= \sum_{k \in \mathbb N_0} \left|\int_{\mathcal O} r(\xi)w(\xi)\chi_{u(k)}(\xi)d\xi\right|^2 =  \int_{\mathcal O} |r(\xi)|^2|w(\xi)|^2 d\xi.
\end{align*}
 Therefore,  we have   condition
\begin{align*}
  \int_{\mathcal O} |r(\xi)|^2|w(\xi)| d\xi= \int_{\mathcal O} |r(\xi)|^2|w(\xi)|^2 d\xi,
\end{align*}
 since  for every $f \in V_\varphi$, we have $\|f\|^2_2 =  \int_{\mathcal O} |r(\xi)|^2|w(\xi)| d\xi. $
That means,
\begin{align*}
  \int_{\mathcal O} |r(\xi)|^2w(\xi)\left(\chi_{\Omega}(\xi)- w(\xi)\right) d\xi=0,
\end{align*}
 holds  for all  integral periodic functions $r \in L^2\left(\mathcal O, w\right)$  if and only if $w(\xi)= \chi_{\Omega}(\xi),  \ a.e.$ \hfill$\square$

Notice that for all $\xi\in K$,   $\displaystyle\sum_{k\in \mathbb N_0}\left|\widehat{\varphi}(\xi+u(k))\right|^2\neq 1,$ if $\widehat{\varphi}=I_{\frak P}$ since $\frak P \subset \mathcal O$ and $\{\mathcal O +u(k):k\in\mathbb N_0\}$ is a measurable partition of $K.$ In the sequel of development of wavelets associated with an MRA on LFPC, Jiang,  Li and  Jin in [18]  found the following result:

\noindent\textbf{{Corollary 3.2.2.}} \textit{A  necessary and sufficient condition to constitute an orthonormal system by $\left\{T_k\varphi: k \in \mathbb N_0\right\}$ is as follows:
\begin{align*}
\sum_{k\in \mathbb N_0}\left|\widehat{\varphi}(\xi+u(k))\right|^2=1, \ \ a.e.,
\end{align*}
for any $\varphi \in L^2(K)$.}

Next, consider a  mapping $\mathcal{T}: L^2(K) \rightarrow L^2\left(\mathcal O,  l^2(\mathcal  Z)\right)$ defined by
$$
\mathcal{T}f(\xi)=\left(\widehat{f}\left(\xi+u(p)\right)\right)_{p \in \mathbb N_0}
$$
which is an isometric isomorphism between $L^2(K) $ and $L^2\left(\mathcal O,  l^2(\mathcal  Z)\right)$ that is an easy consequence of Plancherel theorem. Following is an immediate application of above fiberization map:

 {\noindent\textbf{{Proposition 3.2.3.}}} \textit{Let $V$ be a  $\mathcal Z$-TI subspace of $L^2(K)$. Then the image of $V$ under the fiberization map $\mathcal T$ is  given as follows:
$$
\mathcal{T}(V) =\{ F\in L^2(\mathcal O, l^2(\mathcal Z)): F(\xi)\in J(\xi)\}, \ \mbox{where} \ J(\xi)=\{\mathcal T f(\xi): f \in V\},  \ \mbox{for} \ \xi \in \mathcal O.
$$
The mapping $J$  from $\mathcal O$ to $\{\mbox{closed subspaces of} \ l^ 2(\mathcal  Z)\}$  is known as \textit{range function}, and the dimension  of $J(\xi)$ is  the \textit{multiplicity function} of $V$ which is denoted by $m_V(\xi)(=\mbox{dim} \ J(\xi))$, for $\xi \in \mathcal O$. Throughout, we assume that the range function $J$ is  measurable.}

\textit{Moreover, the multiplicity function satisfies the following condition:
$$
m_{V}(\xi)= \sum_{n=1}^M\sum_{k\in \mathbb N_0} \left|\widehat{\varphi}_n\left(\xi+u(k)\right)\right|^2,
$$
where   $V=\oplus_{i=1}^M V_{\varphi_i}, $ and  $\varphi_i$ is a Parseval frame generator   for the $\mathcal Z$-PTI space $V_{\varphi_i}$.}

\noindent\textbf{{Proof.}} Let $V$ be a  $\mathcal Z$-TI subspace of $L^2(K)$. Then,  there exists a countable set of functions $\{\varphi_i\}_{i=1}^M$ belonging to $V$ such that
$V=\oplus_{i=1}^M V_{\varphi_i},$ where   $M$ is a natural number or infinity. Let us consider a $\mathcal Z$-PTI space $V_{\varphi_i}$ with a Parseval frame generator $\varphi_i$ and a spectrum $\Omega_i$. Then  it is enough to  find the image of   $V_{\varphi_i}$ under the transformation $\mathcal T$. For this,  as every function $f\in V_{\varphi_i}$ satisfies $\widehat{f}(\xi)=r(\xi)\widehat{\varphi_i}(\xi)$, for some integral  periodic function $r \in L^2(\mathcal O, \Omega_i)$,   we have
\begin{align*}
\mathcal{T}(f)(\xi)=&\left(\widehat{f}\left(\xi+u(p)\right)\right)_{p \in \mathbb N_0}=r(\xi)\left(\widehat{\varphi}_i\left(\xi+u(p)\right)\right)_{p \in \mathbb N_0}=r(\xi)\mathcal{T}(\varphi_i)(\xi),
\end{align*}
and hence, we obtain
 \begin{align*}
\mathcal{T}(V_{\varphi_i})=&\left\{F\in L^2\left(\mathcal O, l^2(\mathcal Z)\right): F(\xi)=r(\xi)\mathcal{T}(\varphi_i)(\xi),  r \in L^2\left(\mathcal O, \Omega_i\right)\right\}\\
=& \left\{F\in L^2\left(\mathcal O, l^2(\mathcal  Z)\right): F(\xi)\in J_i(\xi)\right\},
 \end{align*}
where $J_i(\xi)=\mbox{span}\{\mathcal{T}(\varphi_i)(\xi)\}$, {for} $\xi \in \mathcal O$. Therefore, the result follows for $\mathcal Z$-TI space  $V$ by noting  that $J(\xi)=\overline{\mbox{span}}\{\mathcal{T}(\varphi_i)(\xi): i=1, 2, \cdots, M\}=\oplus_{i=1}^M J_{i}(\xi)$.

Moreover,  the multiplicity function satisfies
\begin{align*}
m_{V}(\xi)= \mbox{dim} \  J(\xi)= \sum_{i=1}^M \mbox{dim} \ J_i(\xi)=\sum_{i=1}^M m_{V_{\varphi_{i}}}(\xi)
= \sum_{i=1}^M\sum_{k\in \mathbb N_0} \left|\widehat{\varphi}_i\left(\xi+u(k)\right)\right|^2 ,
\end{align*}
 for a.e. $\xi \in K$. This follows by the  fact that  the dimension of subspace $J_i(\xi)$ of $l^2(\mathcal Z)$  is  one for a.e. $\xi \in \mathcal O \cap \Omega_i$ and  zero  for a.e. $\xi \in \mathcal O$ outside of $\Omega_i$.
\hfill $\square$

In the above theorem, if we apply $\mathcal T^{-1}$ on $\mathcal T(V)$ and using above facts of translation-invariant spaces, we have
$$
V=\{f \in L^2(K): \mathcal  T f(\xi)\in J(\xi), \ \mbox{for  a.e.} \ \xi \in \mathcal O\}.
$$
 It can be easily seen  that the space $V$ is $\mathcal Z$-TI space. 

Analogous to a result of Bownik [7], we state the following:

\noindent{\textbf{Proposition 3.2.4.}} \textit{If $V$ is a $\mathcal Z$-TI space and $N=\|m_V\|_{\infty}$ (N is a natural number or infinity), then there  exist a set of functions $\{\varphi_n\}_{n=1}^N$ such that $V=\oplus_{n=1}^N V_{\varphi_n}$, where $\varphi_n$ is a Parseval frame  generator of $V_{\varphi_n}$.}

Next, we  define spectral function in the context of LFPC which was studied  by Rzeszotnik  in [21] for the case of Euclidean spaces:

\noindent{\textbf{Definition 3.2.5.}} Let $V$ be a $\mathcal Z$-TI space with the range function $J$. Suppose  $P_{J(\xi)}$ is  the orthogonal projection on $J(\xi)$ for a.e. $\xi \in \mathcal O$. Then the \textit{spectral function} of $V$ is the mapping  $\sigma_{V}: K \rightarrow [0, 1]$ given by
$$
\sigma_{V}(\xi+u(k))= ||P_{J(\xi)}e_{u(k)}||_{l^2(\mathcal  Z)}^2, \ \mbox{for} \  \xi~ \in \mathcal O \ \mbox{and} \ k\in \mathbb N_0,
$$
where $\{e_{u(k)}\}_{k\in \mathbb N_0}$ denotes the standard orthonormal basis of $l^2(\mathcal Z)$.



We end this section by summarizing several results of the spectral function by employing above results     whose LCAG version can be found in [11].

\noindent\textbf{Theorem 3.2.6.} \textit{The spectral function satisfies the following properties}:
\begin{enumerate}
	\item [(A)] \textit{Let $V$ be a $\mathcal Z$-TI space of $L^2(K)$ with a set of generators $\{\varphi_n\}^{N}_{n=1}$ ($N$ is  natural number or infinite). If the system $\{T_k \varphi_n: k \in \mathbb N_0, 1 \leq n \leq N\}$  forms a  Parseval frame  for the space $V$, then the function $$\sigma_{V}(\xi)= \sum^N_{n=1} |\widehat{\varphi}_{n}(\xi)|^2,$$ defined for a.e. $\xi \in K$, does not depend on the choice of generators.}
\item [(B)] \textit{Let $\nu$ denote the  set of all $\mathcal Z$-TI spaces of $L^2(K)$. Then there exists a unique mapping $\sigma : \nu \longrightarrow L^\infty(K)$ such that
$$
\sigma_{V_{\varphi}}(\xi)=\begin{cases} |\hat{\varphi}(\xi)|^2\left(\displaystyle{\sum_{k\in \mathbb N_0} |\hat{\varphi}(\xi+u(k))|^2}\right)^{-1}\quad \ &\mbox{for}\quad \xi\in \mbox{supp}\, \hat{\varphi} \\
0   \ &\mbox{otherwise},
\end{cases}
$$
  which is additive on orthogonal sums,
i.e.,    $V=\oplus_{i=1}^{\infty} V_{\varphi_i},$ implies $\displaystyle \sigma_V = \sum^{\infty}_{n=1}\sigma_{V_n}$.}
\item [(C)] \textit{If $V$ is a $\mathcal Z$-TI  space with a decomposition $V=\oplus_{n=1}^N V_{\varphi_n},$ where $\varphi_n$ is a Parseval frame generator for $V_{\varphi_n}$ and $N$ is a natural number or infinity then} $$\sigma_{V}(\xi)= \sum^N_{n=1} |\hat{\varphi}_{n}(\xi)|^2, \  a.e. $$
\item [(D)]\textit{ If  $V$ is a $\mathcal Z$-TI  space,  then  $D_\mathfrak p(V)$ is  $\frak p\mathcal Z$-TI space, and   $\sigma_{D_{\mathfrak p} (V)}(\xi)= \sigma_V({\mathfrak p\xi})$,   a.e.}
\item [(E)] \textit{If  $V$ is a $\mathcal Z$-TI  space,  then $\displaystyle m_V(\xi)=\sum_{k\in \mathbb N_{0}}\sigma_V(\xi+u(k))$,   a.e.}
\end{enumerate}

\noindent\textbf{Proof.} (A)
The result follows  by noting that  $$\sigma_{V}(\xi+u(k))= ||P_{J(\xi)}e_{u(k)}||_{l^2(\mathcal Z)}=\sum^{N}_{n=1}|<P_{J(\xi)}e_{u(k)}, \mathcal T(\varphi_n)(\xi)>_{l^2(\mathcal Z)}|^2$$ $$=\sum^{N}_{n=1}|<e_{u(k)}, \mathcal T(\varphi_n)(\xi)>_{l^2(\mathcal Z)}|^2 =\sum^{N}_{n=1}|\hat{\varphi_n}(\xi+u(k))|^2.$$

\noindent{(B)}   If $V_{\varphi}$ is a $\mathcal Z$-PTI space, then the system $\{T_k \psi\}_{k \in \mathbb N_0}$ is a  Parseval frame for the $\mathcal Z$-PTI space $V_{\varphi}$, where the Fourier transform of $\psi$ is given by
$$
\hat{\psi}(\xi)=\begin{cases} \hat{\varphi}(\xi)\left(\displaystyle {\sum_{k\in \mathbb N_0} |\hat{\varphi}(\xi+u(k))|^2}\right)^{-1/2}\quad \ &\mbox{for}\quad \xi\in \mbox{supp}\, \hat{\varphi} \\
0   \ &\mbox{otherwise},
\end{cases}
$$ and hence,  $\{\mathcal T(\psi)(\xi)\}$ is a Parseval frame for the range function $J_{V_{\varphi}}(\xi)$. So the value of $\sigma_{V_\varphi}$ follows by noting that for all $\xi \in \mathcal O$ and $k \in \mathbb N_0$, we have
$$
\sigma_{V_\varphi}(\xi+u(k))=\|P_{J_{V_{\varphi}}(\xi)}e_{u(k)}\|^2=|<e_{u(k)}, \mathcal T (\psi)(\xi) >|^2=|\hat{\psi}(\xi+u(k))|^2.
$$
The rest portion  of the result follows by considering  $V=\oplus_{n=1}^{\infty} V_{\varphi_n},$ with corresponding range functions $J_V$ and $J_{V_{\varphi_n}}$, and  noting that
   $$\sigma_{V}(\xi+u(k))= ||P_{J_V(\xi)}e_{u(k)}||^{2}_{l^2(\mathcal Z)}=\sum^{\infty}_{n=1}||P_{J_{V_{\varphi_n}}(\xi)}e_{u(k)}||^2
= \sum^{\infty}_{n=1}\sigma_{V_{{\varphi_n}}}(\xi+u(k)). $$

\noindent{(D)} Consider for any $f \in D_\frak p (V)$ and $\gamma \in \frak p \mathcal Z$, i.e., $f=D_{\frak p}g$ for some $g \in V$ and $\frak p^{-1} \gamma=u(l)$, for   $u(l)\in \mathcal Z$. Then for $\xi \in K$, we have
\begin{align*}
f(x-\gamma)= f(x-\frak p u(l))=D_{\frak p}g(x-\frak p u(l))=q^{1/2} g(\frak p^{-1}x-u(l))
=D_{\frak p}(T_lg)(x).
\end{align*}
Since $T_l g \in V$, this shows that $D_{\frak p}(V)$ is $\frak p \mathcal Z$-TI. Next, it suffices to show  the result for $\mathcal Z$-PTI space $V_\varphi$ which has  a  Parseval frame generator $\varphi$. For this, let $f \in D_\frak p (V_\varphi)$. Then, we have
\begin{align*}
||f||^2_2=& \|(D_{\frak p})^{-1}f\|^2=\sum_{k\in \mathbb N_{0}} |<f, D_{\frak p} T_k\varphi>|^2=\sum_{k\in \mathbb N_{0}} \sum_{i=0}^{q-1}|<f, D_{\frak p} T_{qk+i}\varphi>|^2\\
=&\sum_{k\in \mathbb N_{0}}\sum^{q-1}_{i=0}|<f, D_{\frak p}T_{qk}(T_i\varphi)>|^2=\sum_{k\in \mathbb N_{0}}\sum^{q-1}_{i=0}|<f, T_{k}(D_\frak p(T_i\varphi))>|^2,
\end{align*}
in view of the following facts: the map $D_{\frak p}$ is  unitary  on $L^2(K)$; for $k \geq 0$ and $0 \leq i \leq q-1$, $u(qk+i)=\frak p^{-1}u(k)+u(i)$; and $D_{\frak p} T_{qk}=T_{k}D_{\frak p}$, for $k \geq 0$. Since $T_i \varphi \in V_{\varphi}$, this shows  that the system $\{D_{\frak p} (T_i\varphi)\}_{i= 0}^{q-1}$ is a Parseval frame  for the space $D_{\frak p} (V_\varphi)$.

\noindent Note that (C) and (E) follow immediately.
\hfill$\square$

\begin{flushleft}
  \textbf{3.3. Semi-orthogonal Parseval wavelets and GMRAs}
\end{flushleft}

Recall that every semi-orthogonal framelet  $\Psi$  comes from a
GMRA. Now we look into its converse, that is, when a GMRA gives rise to a wavelet, or a semi-orthogonal framelet with the help of knowledge of translation invariant spaces for LFPC obtained from previous subsection.

Following is a main theorem of this section:

\noindent\textbf{Theorem 3.3.1.}  \textit{Suppose that $\Psi$ is a semi-orthogonal Parseval multiwavelet with $L$ generators and $V$ is the space of negative dilates of $\Psi$. Then, $\{D^j_\frak p (V)\}_{j\in\mathbb Z}$ is a GMRA such that $m_{V}(\xi)< \infty $ for a.e. $\xi$, and 
$$
\sum_{d =0}^{q-1} m_{V }\left(\mathfrak p\left(\xi+u(d)\right)\right)- m_{V}(\xi)\leq L, \ \  \mbox{for a.e.}\ \  \xi.
$$
Conversely, if $\{D^j_\frak p (V)\}_{j\in\mathbb Z}$ is a GMRA satisfying above conditions, then there exists  a semi-orthogonal Parseval (multi)wavelet $\Psi$ (with at most $L$ generators) associated with this GMRA.}

   Now,  we proceed as follows to find a proof of above theorem:

\noindent\textbf{Proposition 3.3.2.}  \textit{If $\Psi=\{\psi_l\}_{l=1}^L \subset L^2(K)$ is a semi-orthogonal Parseval multiwavelet with $L$ generators and  $V=\oplus_{j< 0} D^{j}_{\frak p} (W)$, where $W=\overline{\mbox{span}}\{{\psi}_l(\cdot-u(k)):  k \in \mathbb N_0, l=1, 2, \cdots, L\}$,
 then
$$
m_{V}(\xi)=\sum_{l=1}^{L}\sum_{j=1}^{\infty} \sum_{k \in \mathbb N_0}|\widehat{\psi}_l(\mathfrak{p}^{-j} (\xi+u(k)))|^2,  \ \ a.e.,
$$
and hence,  $\int_{\mathcal O} m_{V}(\xi)\leq \frac{L}{q-1}.$ Moreover, it satisfies the consistency equation
$$
m_{V}(\xi)+L \geq \sum_{d =0}^{q-1} m_{V}\left(\mathfrak p\left(\xi+u(d)\right)\right), \ \mbox{for} \ \xi \in K.
$$
}

\noindent\textbf{Proof.}  Since $\displaystyle L^2(K)=V \oplus \oplus_{j\geq 0} D^j_{\frak p}(W)$, we have $\displaystyle \sigma_{V}+\sum_{j\geq 0} \sigma_{D^j_{\frak p}(W)}=1$ in view of Theorem  3.2.6, and hence,  $\displaystyle \sigma_{V }(\xi)=\sum_{l=1}^{L}\sum_{j=1}^{\infty} |\widehat{\psi}_l(\mathfrak{p}^{-j} \xi)|^2.$ This follows by noting  that 
\begin{align*}
\sigma_{V}(\xi)=1-\sum_{j\geq 0} \sigma_{D^j_{\frak p}(W)}(\xi)=\sum_{l=1}^L\sum_{j  \in \mathbb Z} |\hat{\psi}_l(\frak{p}^j \xi)|^2-\sum_{l=1}^L\sum_{j\geq 0} |\hat{\psi}_l(\frak{p}^j \xi)|^2
=\sum_{l=1}^{L}\sum_{j=1}^{\infty} |\widehat{\psi}_l(\mathfrak{p}^{-j} \xi)|^2,
\end{align*}
 because $ \displaystyle\sigma_{D^j_{\frak p}(W)}(\xi)=\sum_{l=1}^L|\hat{\psi}_l(\frak{p}^j \xi)|^2$, for $j \geq 0$.
Therefore, the first result follows by writing    $ \displaystyle m_{V }(\xi)=\sum_{k \in \mathbb N_0} \sigma_{V }(\xi+u(k))$ from Theorem  3.2.6.   Now, we have 
\begin{align*}
\int_{\mathcal O} m_V(\xi) d\xi=&\int_{\mathcal O}  \sum_{l=1}^{L}\sum_{j=1}^{\infty} \sum_{k \in \mathbb N_0}|\widehat{\psi}_l(\mathfrak{p}^{-j} (\xi+u(k))|^2 d\xi
= \sum_{l=1}^{L}\sum_{j=1}^{\infty} \sum_{k \in \mathbb N_0}\int_{\mathcal O+u(k)}  |\widehat{\psi}_l(\mathfrak{p}^{-j} \xi)|^2 d\xi\\
=&\sum_{l=1}^{L}\sum_{j=1}^{\infty}  \frac{1}{q^j}\int_{K}  |\widehat{\psi}_l(\xi)|^2 d\xi
\leq   \frac{L}{q-1}.
\end{align*}
For the consistency equation, we have
\begin{align*}
\sum_{d =0}^{q-1} m_{V}\left(\mathfrak p\left(\xi+u(d)\right)\right)
 =& \sum_{d =0}^{q-1} \sum_{l=1}^{L}\sum_{j=1}^{\infty} \sum_{k \in \mathbb N_0}|\widehat{\psi}_l(\mathfrak{p}^{-j} (\mathfrak p \xi+\mathfrak p u(d)
+u(k)))|^2  \\
=& \sum_{d =0}^{q-1} \sum_{l=1}^{L}\sum_{j=0}^{\infty} \sum_{k \in \mathbb N_0}|\widehat{\psi}_l(\mathfrak{p}^{-j} ( \xi+ u(d)
+\mathfrak p^{-1}u(k)))|^2  \\
=& \sum_{l=1}^{L}  \sum_{k \in \mathbb N_0}|\widehat{\psi}_l( \xi+ u(k))|^2+ \sum_{l=1}^{L} \sum_{j=1}^{\infty} \sum_{k \in \mathbb N_0}|\widehat{\psi}_l(\mathfrak{p}^{-j} ( \xi+ u(k)))|^2  \\
\leq & L+ m_{V}(\xi). \hspace{7cm}  \square
\end{align*}

\noindent\textbf{Corollary 3.3.3.}  \textit{If $f \in L^2(K)$, then the  collection  $\{f(\cdot-u(k)): k\in \mathbb N_0\}$ is a Parseval sequence if and only if $m_f$ satisfies the following consistency equation
$$
\sum_{d =0}^{q-1} m_f\left(\mathfrak p\left(\xi+u(d)\right)\right)\leq 
   1+  m_f(\xi),
$$
where $\displaystyle m_f(\xi)=\sum_{j=1}^{\infty} \sum_{k \in \mathbb N_0}|\widehat{f}(\mathfrak{p}^{-j} (\xi+u(k)))|^2.$}

\noindent\textbf{Corollary 3.3.4.}  \textit{If $f \in L^2(K)$, then the  collection  $\{f(\cdot-u(k)): k\in \mathbb N_0\}$ is an orthonormal system if and only if $m_f$ satisfies the following consistency equation
$$
\sum_{d =0}^{q-1} m_f\left(\mathfrak p\left(\xi+u(d)\right)\right)=
   1+  m_f(\xi),
$$
where $\displaystyle m_f(\xi)=\sum_{j=1}^{\infty} \sum_{k \in \mathbb N_0}|\widehat{f}(\mathfrak{p}^{-j} (\xi+u(k)))|^2.$}

Now, the following result gives a characterization of a   multiwavelet associated with an MRA:

\noindent\textbf{Theorem 3.3.5.} \textit{ If $\Psi=\{\psi_1, \psi_2, \cdots, \psi_{q-1}\} \subset  L^2(\mathbb R)$ is an orthonormal multiwavelet, and $m$ is its associated multiplicity function, then $\Psi$ is   associated with an MRA  if and only if $m\equiv 1$, a.e.}

\noindent\textbf{Proof.}  Suppose $\Psi$ is an MRA multiwavelet. Then   there exists $\varphi \in V$ such that the system  $\{\varphi(\cdot-u(k)):k \in \mathbb N_0\}$ is an orthonormal basis for $V$, and hence,  $\displaystyle \sum_{k \in \mathbb N_0}|\hat{\varphi}(\xi+u(k))|^2=1$, a.e. Therefore,   we have $m(\xi)=1$ since  $V=\overline{\mbox{span}}\{\varphi(\cdot-u(k)):k\in \mathbb N_0\}$ is a  $\mathcal Z$-PTI space.   Conversely, assume that $\Psi \in L^2(\mathbb R)$ is a  multiwavelet, and $m$ is its associated multiplicity function such that $m(\xi)=1$. Then, we have to show that $\Psi$ is an MRA multiwavelet. For this, consider $V=\oplus_{j<0} D^j_{\frak p}(W)$, where   $W=\overline{\mbox{span}}\{{\psi}_l(\cdot-u(k)):  k \in \mathbb N_0, l=1, 2, \cdots, q-1\}$. Then,  the multiplicity function $m(\xi)=1$, and hence,  $V=V_\varphi$, where $\varphi$ is a Parseval frame generator for $V_\varphi$ and spectrum of $V_\varphi$ is equal to $K$ in view of Proposition 3.2.4.  Therefore, we have  $\displaystyle \sum_{k \in \mathbb N_0}|\hat{\varphi}(\xi+u(k))|^2=1$, a.e. and hence,  $\{\varphi(\cdot-u(k)):k\in \mathbb N_0\}$ is an orthonormal basis for $V$.

\hfill $\square$

Next, we define multiplicity function associated with a wavelet as follows:

\noindent{\textbf{Definition 3.3.6.}} If $\Psi$ is a multiwavelet, then there exists a multiplicity function associated to it. This function is called the \textit{wavelet multiplicity function}.

Let $\Psi=\{\psi_l\}_{l=1}^L$ be a  multiwavelet on $L^2(K)$, and consider the wavelet dimension function $D_\Psi(\xi)$ for LFPC  (defined in [5]) which is given by
$$
D_{\Psi}(\xi)=\sum_{l=1}^{L}\sum_{j=1}^{\infty} \sum_{k \in \mathbb N_0}|\widehat{\psi}_l(\mathfrak{p}^{-j} (\xi+u(k)))|^2.
$$

\noindent\textbf{Corollary 3.3.7.}  \textit{Let $\Psi$ be a multiwavelet, and let $m: \mathcal O \rightarrow \mathbb N_0$ be its associated multiplicity function. Then $m(\xi)=D_\Psi(\xi).$}

\noindent\textbf{Proof of Theorem 3.3.1.} Suppose that $\Psi$ is a semi-orthogonal Parseval multiwavelet with $L$ generators and the spaces $V$ and $W$ are defined as follows:\begin{align*} W=\overline{span}\{T_{k}\psi: k\in \mathbb N_{0}, \psi\in \Psi\}, \ \mbox{and} \ V=\overline{span}\{\psi_{j,k}: j<0, k\in \mathbb N_{0}, \psi\in \Psi\}. \end{align*} Then, we have 
\begin{align*}
\int_{\mathcal O}m_{V}(\xi)d\xi=&\int_{K}\sigma_{V}(\xi)d\xi=\sum_{\psi\in \Psi}\sum_{j=1}^{\infty} \int_{K}|\widehat{\psi}(\mathfrak{p}^{-j}\xi|^2 
=\sum_{\psi\in \Psi}||\psi||^2 /(q-1)\\
\leq& L/(q-1)<\infty.
\end{align*}
 Hence, $m_{V}(\xi)< \infty $. Since $W\oplus V= D_\frak p(V)$, We have
$$\sigma_{W}(\xi)+\sigma_{V}(\xi)= \sigma_{D_\frak p(V)}(\xi)=\sigma_{V}(\mathfrak p \xi).$$ This implies that 
$$
m_{W}(\xi)+ m_{V}(\xi)= \sum_{d =0}^{q-1} m_{V }\left(\mathfrak p\left(\xi+u(d)\right)\right).
$$
 Since $m_{W}(\xi)\leq L$, we get the result.

Conversely,  from $$\sum_{d =0}^{q-1} m_{V }\left(\mathfrak p\left(\xi+u(d)\right)\right)- m_{V}(\xi)\leq L$$ and 
 $$m_{W}(\xi)+ m_{V}(\xi)= \sum_{d =0}^{q-1} m_{V }\left(\mathfrak p\left(\xi+u(d)\right)\right),$$ we have $m_{W}(\xi)\leq L$.  By Theorem 3.2.4,  this implies that $W$ has a set $\Psi$ having generators less than or equal to  $L$. Since $V = \bigoplus _{j\leq -1}D^j_\frak p (W)$,
we infer that $\Psi$ is a semi-orthogonal Parseval  multiwavelet associated with the
GMRA $\{D^j_\frak p(V)\}_{j\in\mathbb Z}$. \hfill $\square$

\noindent\textbf{Corollary 3.3.8.}  \textit{Let $V $ be a $\mathcal Z$-TI space such that  the multiplicity function of  $V $ is integrable and satisfies the following consistency equation
$$
m_{V }(\xi)+L=\sum_{d =0}^{q-1} m_{V }\left(\mathfrak p\left(\xi+u(d)\right)\right),
$$
then there exists a set of functions $\Psi=\{\psi_1, \psi_2, \cdots, \psi_L\}$ in $D_{\frak p}(V ) \ominus V \equiv W $ such that the system $\{\psi_l(\cdot-u(k)): k\in \mathbb N_0, 1 \leq l \leq L\}$ is an orthonormal basis for $W $.}

\begin{flushleft}
\textbf{4. Bandlimited Wavelets for LFPC} 
\end{flushleft}

The present section is devoted to the study of  characterizations of  bandlimited   Parseval   multiwavelets as well as Parseval scaling functions for LFPC. 

\noindent\textbf{Proposition  4.1.} \textit{Let $\Psi=\{\psi_m\}_{m=1}^L \subset  L^2(K)$ be such that for each $m\in \{1, 2, \cdots, L\}$, $|\widehat{\psi}_m|=\chi_{W_m},$ and   $W=\bigcup_{m=1}^L W_m$ is a union of measurable subsets of $K$.  Then  $\Psi$ is a  Parseval  multiwavelet in $L^2(K)$ if and only if the following hold}:
\begin{enumerate}
	\item [(i)]  \textit{$\{\mathfrak p^{j} W:  j \in \mathbb Z\}$ is a measurable partition of $K$, and}
		\item [(ii)] \textit{for each $m\in \{1, 2, \cdots, L\}$, the set $\{W_m+u(k):k \in \mathbb N_0\}$ is a measurable partition of   a subset of $K$}.
\end{enumerate}
\textit{In this case, for $j \in \mathbb Z$, $|\frak p^{j}W_m \cap W_{m'}|=0$, where $m, m'\in \{1, 2, \cdots, L\}$, and $m \neq m'$. The set $W$ is known as {Parseval   multiwavelet set of order $L$} in $K$. }

\noindent\textbf{Proof.}    Let $\Psi=\{\psi_m\}_{m=1}^L\subset L^2(K)$ be such that $|\widehat{\psi}_m|=\chi_{W_m},$ where  $W=\bigcup_{m=1}^L W_m$ is a measurable subset of $K$. Then,   the condition (i) of Theorem 3.1.2  yields that $   \bigcup_{j\in\mathbb Z} \mathfrak p^{j} W=K$, a.e.,  that is  equivalent to (i), which also gives that for $j \geq 0$, $|\frak p^{j} W_m \cap W_{m'}|=0$, for each $m, m'\in \{1, 2, \cdots, L\}$, and $m \neq m'$. Further in view of Proposition 3.2.1 and for each $m\in \{1, 2, \cdots, L\}$ ,  the system   $\left\{\psi_m(\cdot-u(k)): k \in \mathbb N_0\right\}$,   is a  Parseval   frame  for $\overline{span}\left\{\psi_m(\cdot-u(k)): k \in \mathbb N_0\right\}$    in $L^2(K)$ if and only if $\sum_{k\in \mathbb N_0}\left|\widehat{\psi}_m (\xi+u(k))\right|^2=\sum_{k\in \mathbb N_0}\chi_{W_m} (\xi+u(k))\leq 1,$     a.e.,  that is  equivalent to the (ii). In this case  $ \{f \in L^2(K): \mbox{supp}\widehat{f} \subset W\}=\overline{\mbox{span}}\{\psi(\cdot-u(k)): \psi \in \Psi, k \in \mathbb N_0\}=:W_0.$  By scaling $W_0$ for any $j \in \mathbb Z$, we have
$
 D^j_\frak p W_0=\overline{\mbox{span}}\{D^j_\frak p \psi(\cdot-u(k)): \psi \in \Psi, k \in \mathbb N_0\}= \{f \in L^2(K): \mbox{supp}\widehat{f} \subset \frak p^{-j} W\}
$. Therefore, $\Psi$ is a Parseval   multiwavelet in $L^2(K)$ if and only if the system $\mathcal A (\Psi)$ forms a Parseval frame for $L^2(K)$ if and only if   $\bigoplus_{ j \in \mathbb Z} D^j_\frak p W_0=L^2(K)$ and (ii) hold, which is true  if and only if (i) and (ii) hold. 
 \hfill $\square$

\noindent\textbf{Corollary 4.2.}  \textit{Let $\Psi=\{\psi_m\}_{m=1}^L   \subset  L^2(K)$ be such that for each $m\in \{1, 2, \cdots, L\}$, $|\widehat{\psi}_m|=\chi_{W_m},$ and   $W=\bigcup_{m=1}^L W_m$ is a union of measurable subsets of $K$.  Then  $\Psi$ is a    multiwavelet in $L^2(K)$ if and only if the following hold}:
\begin{enumerate}
	\item [(i)]  \textit{$\{\mathfrak p^{j} W:  j \in \mathbb Z\}$ is a measurable partition of $K$, and}
		\item [(ii)] \textit{for each $m\in \{1, 2, \cdots, L\}$, the system $\{W_m+u(k):k \in \mathbb N_0\}$ is a measurable partition of     $K$}.		
\end{enumerate}
\textit{In this case, for $j \in \mathbb Z$, $|\frak p^{j}W_m \cap W_{m'}|=0$, where $m, m'\in \{1, 2, \cdots, L\}$, and $m \neq m'$. The set $W$ is known as \textit{ multiwavelet set of order $L$} in $K$.}

 The result given below is a  necessary and sufficient conditions of  Parseval scaling functions for LFPC:

\noindent\textbf{Proposition 4.3.} \textit{Let $\varphi \in L^2(K)$. Then, $\varphi$ is a Parseval scaling function associated to a  PMRA if and only if the following conditions   hold}:
\begin{itemize}
	\item [(i)]\textit{the system $\{ \varphi(\cdot-u(k))\}_{k\in \mathbb N_{0}}$ is a Parseval frame for $\overline{\mbox{span}}\{ \varphi(\cdot-u(k))\}_{k\in \mathbb N_{0}}$ in $L^2(K)$},
	\item [(ii)] \textit{$lim_{j \rightarrow \infty} |\widehat{\varphi}(\mathfrak p^j \xi)|=1$, \ \ a.e. $\xi \in K$, and}
		\item [(iii)] \textit{there exists an integral periodic function $m_0$ in $L^2(\mathcal O)$ such $\widehat{\varphi}(\xi)=m_0(\xi) \widehat{\varphi}(\mathfrak p \xi)$,    a.e. $\xi \in K$}.  
\end{itemize}

\noindent\textbf{Proof.} Those (ii) and (iii) are straightforward   in view of Theorem 5.1 of [4], and (i) follows by noting Proposition 3.2.1. 
\hfill$\square$

Next, we illustrate a characterization of a  Parseval scaling function $\varphi$ such that $|\widehat{\varphi}|=\chi_S,$ for some measurable set $S$ of  $K$. Such set $S$ is known as  \textit{Parseval scaling set}.

\noindent\textbf{Proposition 4.4.} \textit{A function $\varphi$  such that $|\widehat{\varphi}|=\chi_S,$ for some measurable set $S$  of $K$ is a Parseval scaling function of a   PMRA   if and only if} 
\begin{enumerate}
	\item [(i)] $\{S+2k\pi:k\in\mathbb Z^n\}$ \textit{ is a measurable partition of a subset of}  $K$, 	
	\item [(ii)] $\displaystyle \bigcup_{j\in\mathbb Z} {\mathfrak p}^{-j}S=K$, \textit{and} \qquad 	
	 (iii) $S\subset {\mathfrak p}^{-1} S.$
\end{enumerate}
\textit{Moreover,  Parseval  multiwavelet set(s) $W$   associated to  PMRA can be obtained by $W={\mathfrak p}^{-1} S \backslash S$, and hence $S=\displaystyle \bigcup_{j\in\mathbb N} {\mathfrak p}^{j}W$. 
}

\noindent\textbf{Proof.} Suppose $\varphi$  is a Parseval scaling function such that $|\widehat{\varphi}|=\chi_S.$ Then from (i) of Proposition 4.3 and Proposition 3.2.1, we have   
$\sum_{k\in\mathbb N} |\widehat{\varphi} (\xi+u(k))|^2=\chi_{\Omega}(\xi) \Rightarrow \sum_{k\in\mathbb N} \chi_S (\xi+u(k))=\chi_{\Omega}(\xi) \Rightarrow$ $\big(\{S+u(k):k\in\mathbb N\}$, a measurable partition of a subset of $K$.  From  (iii) of Proposition 4.3, $\widehat{\varphi}(\xi)= m(\xi)\widehat{\varphi}(\mathfrak p\xi)\Rightarrow  |\widehat{\varphi}(\xi)|= |m(\xi)||\widehat{\varphi}(\mathfrak p \xi)|\Rightarrow \big(\chi_{S}\leq \chi_{\mathfrak p^{-1}S},$ since $|m(\xi)|\leq1 \big)\Rightarrow  S\subset \mathfrak p^{-1} S$; and  from (ii) of Proposition 4.3, $lim_{j\rightarrow +\infty}|\widehat{\varphi} (\mathfrak p^{j}\xi)|^2=1 \Rightarrow lim_{j\rightarrow +\infty}\chi_S (\mathfrak p^{j}\xi)=1\Rightarrow\big($for every $\epsilon>0,$ there is an $N\in\mathbb N$ such that $|\chi_S (\mathfrak p^{j}\xi)-1|=|\chi_S (\mathfrak p^{j}\xi)-\chi_{K} (\xi)|=|\chi_{(K\backslash \mathfrak p^{-j} S)}|<\epsilon,$  whenever $j>N \big)\Rightarrow \bigcup_{j\in\mathbb Z} \mathfrak p^{-j}S=K$, since $S\subset \mathfrak p^{-1}S$.

Conversely, suppose $\{S+u(k):k\in\mathbb N\}$ is a measurable partition of a subset of  $K$, $\bigcup_{j\in\mathbb Z} \mathfrak p^{-j}S=K$, and $S\subset \mathfrak p^{-1} S.$ Then to prove that $\varphi$  such that $|\widehat{\varphi}|=\chi_S$ is a Parseval scaling function, it suffices to see the conditions (i), (ii) and (iii).  The measurable partition $\{S+u(k):k\in\mathbb N\}$,of a subset of $K$ implies   $\sum_{k\in\mathbb N} \chi_S (\xi+u(k))=\chi_{\Omega}(\xi)\Rightarrow\sum_{k\in\mathbb N} |\widehat{\varphi} (\xi+u(k))|^2=\chi_{\Omega}(\xi);$ $S\subset \mathfrak p^{-1} S \Rightarrow \chi_{S}\leq \chi_{\mathfrak p^{-1}S}\Rightarrow|\widehat{\varphi}(\xi)|= |m(\xi)||\widehat{\varphi}(\mathfrak p\xi)|$, where $|m(\xi)|= \sum_{k\in\mathbb N} \chi_{\mathfrak p S}(\xi+2k\pi)\Rightarrow \widehat{\varphi}(\xi)= m(\xi)\widehat{\varphi}(\mathfrak p\xi)$, where $m(\xi)=\theta(\xi) \sum_{k\in\mathbb N}\chi_{\mathfrak p S}(\xi+u(k))$, for some unimodular, integral periodic function $\theta;$ and $\bigcup_{j\in\mathbb N} \mathfrak p^{-j}S=  K \Rightarrow$ for $j>N, \bigcup_{j>N} \mathfrak p^{-j}S=K,$ since  $S\subset \mathfrak p^{-1} S\Rightarrow lim_{j\rightarrow +\infty}\chi_S (\mathfrak p^{j}\xi)=1\Rightarrow lim_{j\rightarrow +\infty}|\widehat{\varphi} (\mathfrak p^{j}\xi)|^2=1.$ 
\hfill $\square$

\noindent\textbf{Corollary 4.5.} \textit{A function $\varphi$  such that $|\widehat{\varphi}|=\chi_S,$ for some measurable set $S$  of $K$ is a orthonormal scaling function of an   MRA   if and only if }
\begin{enumerate}
	\item [(i)] $\{S+2k\pi:k\in\mathbb Z^n\}$ \textit{is a measurable partition of }   $K$, 	
	\item [(ii)] $\displaystyle \bigcup_{j\in\mathbb Z} {\mathfrak p}^{-j}S=K$, \textit{and} \qquad 	
	 (iii) $S\subset {\mathfrak p}^{-1} S.$
\end{enumerate}
\textit{Moreover,    multiwavelet set(s) $W$ of order $q-1$ associated to  MRA can be obtained by $W={\mathfrak p}^{-1} S \backslash S$, and hence $S=\displaystyle \bigcup_{j\in\mathbb N} {\mathfrak p}^{j}W$, call as orthonormal scaling set. 
}

Next, we provide some examples of orthonormal and Parseval scaling sets.

\noindent\textbf{Example 4.6.} \textbf{(A).} Let us consider Example 3.1.4 in which the multiwavelet set $W$ of order $q-1$ is 
$$W=\displaystyle \bigcup_{j=1}^{q-1} (\mathcal O +u(j))=\mathfrak P^{-1} \backslash \mathcal O.$$ Further, we consider the set $S=\displaystyle \bigcup_{j\in\mathbb N} {\mathfrak p}^{j}W=\displaystyle \bigcup_{j\in\mathbb N} {\mathfrak p}^{j} (\mathfrak P^{-1} \backslash \mathcal O)$. Then, we have $S=\mathcal O$ since $\mathfrak p \mathcal O \subset \mathcal O$ and $\mathfrak P^{-1}=\mathfrak p^{-1} \mathcal O$, and the set $S=\mathcal O$ satisfies conditions (i), (ii) and (iii) of Corollary  4.5. Therefore $S$ is an orthonormal scaling set in the local field  $K$ of positive characteristic, and hence $\Psi$ defined in Example 3.1.4 is associated with an MRA whose scaling function $\varphi$ is defined by $\widehat{\varphi}=\chi_{\mathcal O}$.\\
\textbf{(B).} Let us consider Example 3.1.5  (a) in which the Parseval multiwavelet set $W$ of order $1$ is 
$W=\mathfrak p^m\mathcal O \backslash \mathfrak p^{m+1} \mathcal O,$ where $m \in \mathbb N.$ Further, we consider the set $S=\displaystyle \bigcup_{j\in\mathbb N} {\mathfrak p}^{j}W=\mathfrak p^{m+1} \mathcal O$. Then,   $S$ is a Parseval scaling set in view of  conditions (i), (ii) and (iii) of Proposition  4.4,   and hence $\psi$ defined in Example 3.1.5 (a) is associated with a  PMRA whose Parseval scaling function $\varphi$ is defined by $\widehat{\varphi}=\chi_{\mathfrak p^{m+1} \mathcal O}$.\\
\textbf{(C).} Let us consider Example 3.1.5 (b) in which the Parseval multiwavelet set $W$ of order $q-1$ is 
$W=\displaystyle \bigcup_{j=1}^{q-1} \mathfrak p^m(\mathcal O +u(j))=\mathfrak p^m(\mathfrak P^{-1} \backslash \mathcal O),$ where $m \in \mathbb N.$ Further, we consider the set $S=\displaystyle \bigcup_{j\in\mathbb N} {\mathfrak p}^{j}W=\mathfrak p^{m} \mathcal O$. Then,   $S$ is a Parseval scaling set in view of  conditions (i), (ii) and (iii) of Proposition  4.4,   and hence $\Psi$ defined in Example 3.1.5 (b) is associated with a  PMRA whose Parseval scaling function $\varphi$ is defined by $\widehat{\varphi}=\chi_{\mathfrak p^{m} \mathcal O}$.


\end{document}